\documentclass[12pt,reqno]{amsart}

\begin{document}

\title[Continued fraction expansion]{A continued fraction expansion for a
	$q$-tangent function}

\begin{abstract}

We prove a continued fraction expansion for a certain $q$--tangent
function that was conjectured by Prodinger.

\end{abstract}

\author{Markus Fulmek}
\address{Institut f\"ur Mathematik der Universit\"at Wien\\
Strudlhofgasse 4, A-1090 Wien, Austria}
\email{{\tt Markus.Fulmek@Univie.Ac.At}\newline\leavevmode\indent
{\it WWW}: {\tt http://www.mat.univie.ac.at/\~{}mfulmek}
}
\date{\today}

\rm

\maketitle

\def\qp#1#2{[{#1}^{#2}]}
\def\sqp#1#2{s_{\qp{#1}{#2}}}
\def\pa#1{[#1]}
\def\sla#1{s_{\pa{#1}}}
\def\TP{P}
\def\TQ{\hat{P}}
\def\TS{\check{P}}
\def\minor#1#2#3#4#5{#1_{\{#2,#3\},\{#4,#5\}}}
\def\la{\lambda}
\def\si{\sigma}
\def\pprime{{\prime\prime}}

\def\q#1{[#1]_q}
\def\pq#1#2#3{(#1;#3)_{#2}}
\def\ceil#1{\lceil #1\rceil}
\def\floor#1{\lfloor #1\rfloor}
\newcommand{\HypFsimple}[2]{\sideset{_#1}{_#2}{F}}
\newcommand{\HypPhi}[6]{\sideset{_#1}{_#2}{\Phi}\!\left[\matrix  
#3\\#4\endmatrix;{\displaystyle #5,#6}\right]}

\let\epsilon\varepsilon
\newcommand{\esf}{{e}}                  
\newcommand{\hsf}{{h}}                  
\newcommand{\ksf}{{k}}                  
\newcommand{\schf}[2]{{s}_{#1}\of{#2,\bold{x}}}
\newcommand{\defeq}{:=}                                 
\newcommand{\tableau}{T}                                
\newcommand{\Ferrer}{B}                            
\newcommand{\digraph}{D}                                
\newcommand{\sdigraph}{S}                               
\newcommand{\horedge}{a}                                
\newcommand{\veredge}{\hat{a}}          
\newcommand{\path}{P}                                   
\newcommand{\pathset}{M}                                
\newcommand{\fpath}{{\mathcal P}}                       
\newcommand{\nfpath}{{\cal N}}          
\newcommand{\nifpath}{{\cal P}_o}       
\newcommand{\SG}[1]{\boldsymbol{S}_{#1}}
\newcommand{\GF}[1]{\bold{GF}\of{#1}}
\newcommand{\sgn}{\operatorname{sgn}}                   
\newcommand{\inv}{\operatorname{inv}}                   
\newcommand{\sschf}[2]{{sp}_{#1}\of{#2,\bold{x}}}
\newcommand{\isschf}[2]{{sp}_{n,m}\of{#1,#2;\bold{x},\bold z}}
\newcommand{\osschf}[2]{{sp}_{n,1}\of{#1,#2,\bold{x},z}}
\newcommand{\reflect}{{\bold{R}}}               
\newcommand{\mreflect}{\tilde{\bold{R}}}                
\newcommand{\oschf}[2]{{o}_{#1}\of{#2,\bold{x}}}
\newcommand{\oschftriv}[2]{{o}_{#1}\of{#2,\seqof{1,1,\dots,1}}}
\newcommand{\liff}{\text{ iff }}                
\newcommand{\limp}{\text{ implies }}
\renewcommand{\land}{\text{ and }}      

\newcommand{\thmref}[1]{Theorem~\ref{#1}}
\newcommand{\secref}[1]{\S\ref{#1}}
\newcommand{\lemref}[1]{Lemma~\ref{#1}}
\newcommand{\figref}[1]{Figure~\ref{#1}}

\newcommand{\seqof}[1]{\left\langle#1\right\rangle}
\newcommand{\setof}[1]{\left\{#1\right\}}
\newcommand{\of}[1]{\left(#1\right)}
\newcommand{\parof}[1]{\left(#1\right)}
\newcommand{\numof}[1]{\left|#1\right|}
\newcommand{\detof}[2]{\left|{#2}\right|_{#1\times #1}}
\newcommand{\firstcolumn}[1]{
\quad \vdots\quad }
\newcommand{\firstrow}[1]{
\quad \vdots\quad }
\newcommand{\brkof}[1]{\left[#1\right]}
\newcommand{\intof}[1]{\lfloor #1\rfloor}
\newcommand{\absof}[1]{\left|#1\right|}
\newcommand{\sgnof}[1]{\operatorname{sgn}\of{#1}}

\newtheorem{thm}{Theorem}
\newtheorem{lem}[thm]{Lemma}
\newtheorem{dfn}[thm]{Definition}
\newtheorem{obs}[thm]{Observation}
\newtheorem{rem}[thm]{Remark}


\section{Introduction}
\label{intro}

In \cite{Prodinger}, Prodinger defined the following $q$--trigonometric functions
\begin{align*}
\sin_q(z)&=\sum_{n=0}^\infty\frac{(-1)^nz^{2n+1}}{\q{2n+1}!}q^{n^2},
	\\
\cos_q(z)&=\sum_{n=0}^\infty\frac{(-1)^nz^{2n}}{\q{2n}!}q^{n^2}.
\end{align*}
Here, we use standard $q$--notation:
\begin{gather*}
\q{n}:=\frac{1-q^n}{1-q},\;\q{n}!:=\q{1}\q{2}\dots\q{n},\\
	\pq{a}{n}{q}:=(1-a)(1-a q)\dots(1-a q^{n-1}).
\end{gather*}
These $q$--functions are variations of Jackson's \cite{Jackson}
$q$-sine and $q$-cosine functions.

For the $q$--tangent function $\tan_q=\frac{\sin_q}{\cos_q}$, Prodinger
conjectured the following continued fraction expansion
(see \cite[Conjecture 10]{Prodinger}):
\begin{equation}
\label{eq:prodinger}
-z\tan_q(z)
=
-\cfrac{z^2}{\q{1}q^{0}-
	\cfrac{z^2}{\q{3}q^{-2}-
		\cfrac{z^2}{\q{5}q^{1}-
			\cfrac{z^2}{\q{7}q^{-9}-\dotsb
			}
		}
	}
}.
\end{equation}
Here, the powers of $q$ are of the form ${(-1)^{n-1} n(n-1)/2-n+1}$.

The purpose of this note is to prove this statement. In our proof, we make
use of the polynomials (see \cite[\S2, (11)]{Perron}) $A_n(z)$ and $B_n(z)$,
which are given recursively by
\begin{align}
A_n(z)&=b_n A_{n-1}(z) -z^2A_{n-2}(z)\label{eq:Arec},\\
B_n(z)&=b_n B_{n-1}(z) -z^2B_{n-2}(z)\label{eq:Brec};
\end{align}
with initial conditions (see \cite[\S2, (12)]{Perron})
\begin{gather*}
A_{-1}=1,\;B_{-1}=0,\;A_0=b_0,\;B_0=1,
\end{gather*}
where $b_0=0$, $b_n=\q{2n-1}q^{(-1)^{n-1} n(n-1)/2-n+1}$. As is well known (see
\cite[\S2]{Perron}), the continued fraction terminated after the term $b_n$
is equal to $\frac{A_n}{B_n}$, whence \eqref{eq:prodinger} follows from the
assertion
\begin{equation}
\label{eq:continuants}
A_n\cos_q+z B_n\sin_q=O(z^{2n+1}),
\end{equation}
i.e., the leading $2n$ coefficients of $z$ vanish in \eqref{eq:continuants}.

In Section~\ref{sec:proof} we give a proof of \eqref{eq:continuants}
(and thus of \eqref{eq:prodinger}).


\section{The proof}

\label{sec:proof}
Both $A_n$ and $B_n$ are polynomials in $z^2$:
\begin{gather*}
A_n(z)=\sum_j c_{n,j}z^{2j},\;
B_n(z)=\sum_j d_{n,j}z^{2j}.
\end{gather*}
Observe that from the recursions \eqref{eq:Arec} and \eqref{eq:Brec} we
obtain immediately
\begin{gather}
c_{n,k}=b_n c_{n-1,k}-c_{n-2,k-1}\text{ and }
d_{n,k}=b_n d_{n-1,k}-d_{n-2,k-1},\label{eq:cdrec}
\end{gather}
with initial conditions
\begin{gather*}
c_{0,k}=d_{0,k}=c_{-1,k}=d_{-1,k}=c_{n,0}=d_{n,-1}=0,\text{ }
c_{1,1}=-1,\text{ }d_{0,0}=1.
\end{gather*}

Given this notation, we have to prove the following assertion for the
coefficients of $z^{2k}$ in \eqref{eq:continuants}: For $n\geq 1$,
$0\leq k\leq n$, there holds
\begin{equation}
\label{eq:toproof}
\sum_{i=0}^k c_{n,i}\frac{(-1)^{k-i}}{\q{2k-2i}!}q^{(k-i)^2} + 
\sum_{i=0}^k d_{n,i}\frac{(-1)^{k-i-1}}{\q{2k-2i-1}!}q^{(k-i-1)^2} = 0.
\end{equation}
In fact, we shall state and prove a slightly more general assertion:

\begin{lem}
\label{lem:Lemma}
Given the above definitions, we have for all $n\geq 1$, $k\geq 0$:
\begin{multline}
\label{eq:toproof2}
\sum_{i=0}^{k-1}(-1)^i\frac{q^{(k-i-1)^2}}{\q{2k-2i-2}!}\left(
	c_{n,i+1}+\frac{d_{n,i}}{\q{2k-2i-1}}
\right) = \\
(-1)^n q^{(5+3(-1)^n-12k-4(-1)^n k+8k^2+8n-8k n+4n^2-2(-1)^n n^2)/8}\\
\times\frac{
	\prod_{s=k-n}^k\q{2s}
}{
	\q{2k}!
}.
\end{multline}
\end{lem}
Note that the left hand side of \eqref{eq:toproof2} is the same as in
\eqref{eq:toproof}, and the right hand side of \eqref{eq:toproof2} vanishes
for $0\leq k\leq n$. Hence \eqref{eq:toproof} (and thus Prodinger's conjecture)
is an immediate consequence of Lemma~\ref{lem:Lemma}. 

\begin{proof}
We perform an induction on $k$ for arbitrary $n$.

The case $k=0$ is immediate. For the case $k=1$, observe that
\begin{gather*}
-c_{n,1} = d_{n,0} = \prod_{s=1}^n b_n.
\end{gather*}

For the inductive step $(k-1)\rightarrow k$, we shall rewrite the recursions
\eqref{eq:cdrec} in the following way:
\begin{gather*}
c_{n,k}=-\sum_{i=0}^{n-2}\left(c_{i,k-1}\prod_{j=i+3}^n b_j\right),\;
d_{n,k}=-\sum_{i=0}^{n-2}\left(d_{i,k-1}\prod_{j=i+3}^n b_j\right).
\end{gather*}

Substitution of these recursions into \eqref{eq:toproof2} and interchange of
summations transform the identity into
\begin{equation*}
\frac{
q^{(k-1)^2}(1-\q{2k-1})\prod_{s=1}^n b_s
}{
\q{2k-1}!
} +
\sum_{i=0}^{n-2}\left({\text{rhs}}(i,k-1)\prod_{j=i+3}^n b_j\right)  =
	{\text{rhs}}(n,k),
\end{equation*}
where $\text{rhs}(n,k)$ denotes the right hand side of \eqref{eq:toproof2}.

Now we use the induction hypothesis. As it turns out, factorization of powers
of $q$ from $\left({\text{rhs}}(i,k-1)\prod_{j=i+3}^n b_j\right)$ yields
the same power
for $2i$ and $2i+1$, whence we can group these terms together. After several
steps of simplification we arrive at the following identity:
\begin{multline}
\label{eq:identity}
\left(
\sum_{j=0}^{\ceil{\frac{n-4}{2}}}\frac{
	\pq{q^{-2k+6}}{j}{q^4}\pq{q^{-2k+4}}{j}{q^4}
		\pq{q^{17/2-k}}{j}{q^4}\pq{-q^{17/2-k}}{j}{q^4}
	}{
	\pq{q^7}{j}{q^4}\pq{q^9}{j}{q^4}\pq{q^{9/2-k}}{j}{q^4}\pq{-q^{9/2-k}}{j}{q^4}
	}q^{(2k-1)j}
\right)
\\
\times\frac{
	\pq{q}{n}{q^2}q(1-q^{2k-1})(1-q^{2k-2})(1-q^{9-2k})
}{
	(1-q)(1-q^3)(1-q^5)
}
-(1-q^{2k-1})\pq{q^3}{n-1}{q^2}\\
-(-1)^nq^{(-1+(-1)^n+2k-2(-1)^nk+2n-4kn+4n^2)/4}\pq{q^{2k-2n}}{n}{q^2}\\
+\pq{q}{n}{q^2}+
	\chi(n)(1-q^{2k-1})q^{n(2n-2k+1)/2}\pq{q^{2k-2n+2}}{n-1}{q^2}
= 0,
\end{multline}
where $\chi(n)=1$ for $n$ even and $0$ for $n$ odd.

The sum can be evaluated by means of the very--well--poised
$\sideset{_6}{_5}{\operatorname{\phi}}$ summation formula
\cite[(2.7.1); Appendix (II.20)]{Gasper}:
\begin{multline}
\label{eq:hypsummation}
\sum_{j=0}^{\infty}
\frac{\pq{a}{j}{q}\pq{\sqrt{a} q}{j}{q}\pq{-\sqrt{a} q}{j}{q}
	\pq{b}{j}{q}\pq{c}{j}{q}\pq{d}{j}{q}}{
\pq{q}{j}{q}\pq{\sqrt{a}}{j}{q}\pq{-\sqrt{a}}{j}{q}
	\pq{\frac{a q}{b}}{j}{q}\pq{\frac{a q}{c}}{j}{q}\pq{\frac{a q}{d}}{j}{q}}
\left(\frac{a q}{b c d}\right)^j
\\=
	\frac{
	\pq{a q}{\infty}{q}\pq{\frac{a q}{b c}}{\infty}{q}
	\pq{\frac{a q}{b d}}{\infty}{q}\pq{\frac{a q}{c d}}{\infty}{q}
	}{
	\pq{\frac{a q}{b}}{\infty}{q}\pq{\frac{a q}{c}}{\infty}{q}
	\pq{\frac{a q}{d}}{\infty}{q}\pq{\frac{a q}{b c d}}{\infty}{q}
	}.
\end{multline}
The sum we are actually interested in does not extend to infinity, so we
rewrite is as follows:
\begin{align*}
\sum_{j=0}^{\ceil{\frac{n-4}{2}}} s(n,k,j)
&=\sum_{j=0}^{\infty} s(n,k,j) -
\sum_{j=\ceil{\frac{n-2}{2}}}^{\infty} s(n,k,j)\\
&=\sum_{j=0}^{\infty} s(n,k,j) -
s(n,k,\ceil{\frac{n-2}{2}})\sum_{j=0}^{\infty}
	\frac{s(n,k,j+\frac{n-2}{2})}{s(n,k,\ceil{\frac{n-2}{2}})},
\end{align*} 
where $s(n,k,j)$ denotes the summand in \eqref{eq:identity}.
Now, replacing $q$ by $q^4$, $a$ by $q^{-2k+8a+9}$, $b$ by $q^{-2k+4a+4}$,
$c$ by $q^{-2k+4a+6}$ and $d$ by $q^{4}$ in the summand of
\eqref{eq:hypsummation} gives $\frac{s(n,k,j+a)}{s(n,k,a)}$ times the fraction $\frac{\pq{q^{-2k+8a+9}}{j}{q^4}\pq{q^4}{j}{q^4}}
{\pq{q^{-2k+8a+9}}{j}{q^4}\pq{q^4}{j}{q^4}}$, which cancels.
So we obtain after some simplification:
\begin{multline*}
\sum_{j=0}^{x}\frac{
	\pq{q^{-2k+6}}{j}{q^4}\pq{q^{-2k+4}}{j}{q^4}
		\pq{q^{17/2-k}}{j}{q^4}\pq{-q^{17/2-k}}{j}{q^4}
	}{
	\pq{q^7}{j}{q^4}\pq{q^9}{j}{q^4}\pq{q^{9/2-k}}{j}{q^4}\pq{-q^{9/2-k}}{j}{q^4}
	}q^{(2k-1)j}\\=
\frac{(1-q^{3})(1-q^{5})}{(1-q^{2k-1})(1-q^{9-2k})}
\left
	(1-q^{(2k-1)(x+1)}\frac{\pq{q^{-2k+4}}{2x+2}{q^2}}{\pq{q^3}{2x+2}{q^2}}
\right).
\end{multline*}

Substitution of this evaluation in \eqref{eq:identity} and simplification
yield for
both cases $n$ even ($n=2N$) and odd ($n=2N-1$) the same equation
\begin{equation*}
\pq{q^{-2k+2}}{2N-1}{q^2}=-q^{-2(k-N)(2N-1)}\pq{q^{2k-4N+2}}{2N-1}{q^2},
\end{equation*}
which, of course, is true. This finishes the proof.
\end{proof}


\ifx\undefined\bysame
\newcommand{\bysame}{\leavevmode\hbox to3em{\hrulefill}\,}
\fi


\end{document}